\documentclass[12pt, leqno,twoside]{article}
\usepackage{amssymb}
\usepackage{amsmath}
\usepackage{amsthm}
\usepackage{graphicx}
\usepackage[pagebackref,colorlinks,citecolor=blue,linkcolor=blue]{hyperref}
\usepackage{cite}   

\usepackage{ulem}
\usepackage{amsmath}
\usepackage{amsfonts}
\usepackage{amssymb}
\usepackage{amsmath,bbm,amssymb,amsxtra}
\usepackage{mathrsfs}
\usepackage{enumerate}
\usepackage{caption}
\allowdisplaybreaks
\usepackage{epsfig}
\usepackage{ae}
\def\XXint#1#2#3{{\setbox0=\hbox{$#1{#2#3}{\int}$ }
		\vcenter{\hbox{$#2#3$ }}\kern-.6\wd0}}

\newcommand{\loc}{\text{loc}}

\newcommand{\mbd}{{\mathbb D}}
\newcommand{\mbx}{{\mathbb X}}

\newcommand{\mbc}{{\mathbb C}}
\newcommand{\mby}{{\mathbb Y}}
\newcommand{\mbs}{{\mathbb S}}
\newcommand{\mbr}{{\mathbb R}}

\newcommand{\X}{{\mathbb X}}
\newcommand{\pao}{{\partial \Omega}}

\newcommand{\D}{{\mathbb D}}

\newcommand{\real}{{\mathbb R}}

\newcommand{\rarrow}{\rightarrow}
\newcommand{\ra}{\rightarrow}
\newtheorem{thm}{Theorem}[section]
\newtheorem{lem}[thm]{Lemma}

\newtheorem{cor}[thm]{Corollary}
\newtheorem{defn}[thm]{Definition}
\newtheorem{example}[thm]{Example}
\numberwithin{equation}{section}
\newtheorem*{thm1}{Theorem A}

\makeatletter
\newcommand*{\mint}[1]{%
  \mint@l{#1}{}%
}
\newcommand*{\mint@l}[2]{%
  \@ifnextchar\limits{%
    \mint@l{#1}%
  }{%
    \@ifnextchar\nolimits{%
      \mint@l{#1}%
    }{%
      \@ifnextchar\displaylimits{%
        \mint@l{#1}%
      }{%
        \mint@s{#2}{#1}%
      }%
    }%
  }%
}
\newcommand*{\mint@s}[2]{%
  \@ifnextchar_{%
    \mint@sub{#1}{#2}%
  }{%
    \@ifnextchar^{%
      \mint@sup{#1}{#2}%
    }{%
      \mint@{#1}{#2}{}{}%
    }%
  }%
}
\def\mint@sub#1#2_#3{%
  \@ifnextchar^{%
    \mint@sub@sup{#1}{#2}{#3}%
  }{%
    \mint@{#1}{#2}{#3}{}%
  }%
}
\def\mint@sup#1#2^#3{%
  \@ifnextchar_{%
    \mint@sup@sub{#1}{#2}{#3}%
  }{%
    \mint@{#1}{#2}{}{#3}%
  }%
}
\def\mint@sub@sup#1#2#3^#4{%
  \mint@{#1}{#2}{#3}{#4}%
}
\def\mint@sup@sub#1#2#3_#4{%
  \mint@{#1}{#2}{#4}{#3}%
}
\newcommand*{\mint@}[4]{%
  \mathop{}%
  \mkern-\thinmuskip
  \mathchoice{%
    \mint@@{#1}{#2}{#3}{#4}%
        \displaystyle\textstyle\scriptstyle
  }{%
    \mint@@{#1}{#2}{#3}{#4}%
        \textstyle\scriptstyle\scriptstyle
  }{%
    \mint@@{#1}{#2}{#3}{#4}%
        \scriptstyle\scriptscriptstyle\scriptscriptstyle
  }{%
    \mint@@{#1}{#2}{#3}{#4}%
        \scriptscriptstyle\scriptscriptstyle\scriptscriptstyle
  }%
  \mkern-\thinmuskip
  \int#1%
  \ifx\\#3\\\else_{#3}\fi
  \ifx\\#4\\\else^{#4}\fi
}
\newcommand*{\mint@@}[7]{%
  \begingroup
    \sbox0{$#5\int\m@th$}%
    \sbox2{$#5\int_{}\m@th$}%
    \dimen2=\wd0 %
    \let\mint@limits=#1\relax
    \ifx\mint@limits\relax
      \sbox4{$#5\int_{\kern1sp}^{\kern1sp}\m@th$}%
      \ifdim\wd4>\wd2 %
        \let\mint@limits=\nolimits
      \else
        \let\mint@limits=\limits
      \fi
    \fi
    \ifx\mint@limits\displaylimits
      \ifx#5\displaystyle
        \let\mint@limits=\limits
      \fi
    \fi
    \ifx\mint@limits\limits
      \sbox0{$#7#3\m@th$}%
      \sbox2{$#7#4\m@th$}%
      \ifdim\wd0>\dimen2 %
        \dimen2=\wd0 %
      \fi
      \ifdim\wd2>\dimen2 %
        \dimen2=\wd2 %
      \fi
    \fi
    \rlap{%
      $#5%
        \vcenter{%
          \hbox to\dimen2{%
            \hss
            $#6{#2}\m@th$%
            \hss
          }%
        }%
      $%
    }%
  \endgroup
}

\begin{document}
\title{\Large\bf  Improved regularity of harmonic diffeomorphic extensions on quasihyperbolic domains
\footnotetext{\hspace{-0.35cm}
 $2010$ {\it Mathematics Subject classfication}: 58E20, 46E35, 30C62
 \endgraf{{\it Key words and phases}: Poisson extension, Orlicz-Sobolev homeomorphisms, weighted Sobolev homeomorphisms, quasihyperbolic domains }
\endgraf{{\it ${}^{\mathbf{*}}$ Corresponding author}}
 }
}
\author{Zhuang Wang, Haiqing Xu${}^{\mathbf{*}}$}
\date{ }
\maketitle

\begin{abstract}
Let $\mbx$ be a Jordan domain satisfying hyperbolic growth conditions.
Assume that $\varphi$ is a homeomorphism from the boundary $\partial \mbx$ of $\mbx$ onto the unit circle. Denote by $h$ the harmonic diffeomorphic extension of $\varphi $ from $\mbx$ onto the unit disk. We establish the optimal Orlicz-Sobolev regularity and weighted Sobolev estimate of $h.$ These generalize the Sobolev regularity of $h$ by Koski-Onninen \cite[Theorem 3.1]{Koski 2018}. 
\end{abstract}

\section{Introduction}

A planar Jordan curve $\Gamma$ is a simple closed curve, i.e. a non-self-intersection continuous curve, in the plane $\mbr^2$. The famous \textit{Jordan Curve Theorem} asserts that each planar Jordan curve $\Gamma$ divides the plane into two regions (an ``interior'' and an ``exterior''), with $\Gamma$ being the common boundary. Moreover every continuous path, connecting a point in one region to a point in the other, intersects with $\Gamma$ somewhere.
Every region with a Jordan curve as its boundary is called a \textit{Jordan domain}.
Due to its importance in low-dimensional topology and complex analysis, the Jordan Curve Theorem received great attention from prominent mathematicians of the first half of the 20th century. Various proofs of the theorem and its generalizations were obtained by J. W. Alexander, L. Antoine, L. Brouwer, and A. Pringsheim et al. One generalization from A. Schoenflies is as follows: 
\begin{thm1}
Fix two Jordan domains $\mbx \subset \mbr^2$ and $\mby \subset \mbr^2 .$
Every homeomorphism $\varphi$ from the boundary $\partial \mbx$ of $\mbx$ onto $\partial \mby$ has a homeomorphic extension $\Phi :\mbr^2 \ra \mbr^2 .$
\end{thm1}
We shall refer to Theorem A as the \textit{Jordan-Schoenflies theorem}.
In geometric topology, people concern diffeomorphic approximation of $\Phi$, see \cite{Moise 1997}.
On the other hand, it is more important for analysts to explore the regularity of $\Phi .$  
By \cite[Example 4.3]{Xu 2020}, there is a self-homeomorphism $\varphi$ on the unit circle $\mbs .$ From its construction, we find that $\varphi$ fails to be differentiable on a subset of $\mbs$. Therefore $\varphi$ does not have a diffeomorphic extension. Meanwhile note that $\varphi$ has weak derivative on the whole $\mbs.$ One can expect that its extension $\Phi$ is weakly differentiable on the entire $\mbr^2.$ Hence it is reasonable to discuss the Sobolev $W^{1,p}$ regularity of $\Phi.$

Fix Jordan domains $\mbx \subset \mbr^2$ and $\mby \subset \mbr^2 .$ Let
\begin{equation*}
H^{1,p} (\mbx ,\mby)= \left\{h : h: \mbx \ra \mby \mbox{ is a homeomorphism in } W^{1,p} (\mbx ,\mby) \right\}.
\end{equation*}
J. M. Ball\cite{Ba76, Ball 1981} introduced $H^{1,p} (\mbx ,\mby)$ to study mathematical problems of nonlinear elasticity. Denote by $f$ a deformation of a body in $\mbr^n .$ 
\begin{enumerate}
\item[(1)]During a deformation in nonlinear elasticity, the material cannot break and no cavities are created. This corresponds to that $f$ is continuous.
\item[(2)]No interpenetration of material and invertibility of deformation process require that $f$ is injective and invertible, respectively.
\item[(3)]Not all deformation $f$ is a diffeomorphism. Therefore when we study energy of a deformation homeomorphism $f$, the Sobolev regularity of $f$ is a natural choice. 
\end{enumerate}

After Ball's work, people begin to use calculus of variations to study mathematical problems of nonlinear elasticity. We summarize these problems as follows: find the minimization of an energy functional for a class of Sobolev deformation homeomorphisms, i.e.
\begin{equation}\label{0323}
\inf_{h \in H^{1,p} (\mbx ,\mby)} \int_{\mbx} E(x,h(x),Dh (x)) \, dx .
\end{equation}  
From the physical background, physicists often directly assume that $H^{1,p} (\mbx ,\mby)$ in \eqref{0323} is not empty.
For mathematicians, it is interesting to offer a rigorous mathematical proof.    
By a boundary homeomorphism we can describe a deformation from $\mbx$ onto $\mby .$ Therefore to check that $H^{1,p} (\mbx ,\mby) \neq \emptyset$, it suffices to show that for a homeomorphism $\varphi :\partial \mbx \ra \partial \mby$ as in the Jordan-Schoenflies theorem, there is a Sobolev homeomorphic extension. We call such a problem the \textit{Sobolev-Jordan-Schoenflies problem}.

Function theory of $H^{1,p} (\mbx, \mby)$ is also interesting.  
The Ball-Evans problem asks for the diffeomorphic approximation of $H^{1,p} (\mbx, \mby) $ mappings. 
Its solution will simplify the function classes in the minimization problem \eqref{0323}.
See \cite{Iwaniec 2011, Hencl 2018, Campbell 2018} for answers to the Ball-Evans problem. 
Results on the Sobolev-Jordan-Schoenflies problem belong to the trace theory of $H^{1,p} (\mbx, \mby) .$
The trace theory of classical Sobolev spaces has important applications in the boundary value problem of PDE. We expect that the trace theory of $H^{1,p} (\mbx, \mby) $ can be potentially applied to analogous problems.

Let $\mby$ be a bounded convex Jordan domain in the plane. Assume that $\varphi$ is a homeomorphism from the unit circle $\mbs$ onto the boundary $\partial \mby $ of $\mby .$ Denote by $h$ the complex-valued Poisson extension of $\varphi$, that is 
\begin{equation*}
h(z) = \frac{1}{2 \pi}\int_{\mbs} \frac{1-|z|^2}{|z-\xi|^2} \varphi (\xi) \, |d \xi| 
\end{equation*}
for all $z \in \mbd .$
Then the Rad\'o-Kneser-Choquet theorem \cite{Rado 1926, Kneser 1926, Choquet 1945} states that $h$ is a univalent mapping from $\mbd$ onto $\mby .$ 
In view of the ``averaging'' property of $h,$ the requirement that $\mby$ is convex is necessary.
Furthermore Lewy showed in \cite{Lewy 1936} that the Jacobian of $h$ nowhere vanishes. Hence $h$ is a real analytic diffeomorphism from $\mbd$ onto $\mby .$ 
We refer to \cite{Duren 2004} written by P. Duren for more on complex-valued harmonic mappings in the plane.
Notice that derivatives of $h$ are not always uniformly bounded. Therefore when we study the regularity of derivatives of $h$ on $\mbd,$ the Sobolev regularity of $h$ becomes the theme.

The study for Sobolev regularity of the harmonic extension $h$ stems from Astala et al. \cite{Astala 2005}. With additionally the $C^1$-smoothness of $\partial \mby,$ they showed that the mapping from $\mby$ onto $\mbd$, whose distortion function has the minimal $L^1$ integrability, is exactly the inverse of the unique harmonic extension $h : \mbd \ra \mby .$ Furthermore this minimal $L^1$ energy equals to the square integrability of the derivative of $h$ on $\mbd,$ which is comparable to that
\begin{equation}\label{1230}
\int_{\partial \mby} \int_{\partial \mby} |\log|\varphi^{-1} (\xi) -\varphi^{-1} (\eta)|| \, |d \xi| \, |d \eta|< \infty .
\end{equation} 
The condition \eqref{1230} does not automatically hold. In fact Verchota \cite{Verchota 2007} showed that the harmonic homeomorphisms of $\mbd$ onto itself need to be in the Sobolev space $W^{1,p} (\mbd, \mbd)$ for $p<2,$ but not in $W^{1,2} .$ Iwaniec et al. \cite{Iwaniec 2009} established more delicate estimates, which provided a weak type $L^2$-estimate and estimates in Orlicz classes near $L^2 (\mbd) $ for the gradient of $h .$ 
Later Xu et al. \cite{our paper, Xu 2020} generalized these results to the case that $\mby$ is an internal chord-arc Jordan domain.
Their extension mapping $\widetilde h :\mbd \ra \mby$ has the form
\begin{equation*}
\widetilde h= h \circ F,
\end{equation*}
where $h :\mbd \ra \mbd$ is the Poisson extension and $F :\mbd \ra \mby$ is a bi-Lipschitz mapping with respect to the internal metric $\lambda_{\mby}.$ Here $\lambda_{\mby} (y_1,y_2)$ for $y_1 \in \mby$ and $y_2 \in \mby$ is defined as the infimum of lengths of all rectifiable curves in $\mby$ joining $y_1$ to $y_2 .$ 
Moreover they studied connections between $W^{1,p} (\mbd, \mby)$ regularity of $\widetilde h$ for $p \ge 1,$ a double integral condition on $\varphi^{-1}$ like \eqref{1230}, and the internal $p$-Douglas condition
\begin{equation*}
\int_{\mbs} \int_{\mbs} \frac{(\lambda_{\mby} (\varphi (\xi),\varphi(\eta)))^p}{|\xi -\eta|^p} \, |d \xi|\, |d \eta|<\infty .
\end{equation*}
They also explored weighted Sobolev regularity and Orlicz-Sobolev regularity of the extension mapping $\widetilde h$.

In all above works, harmonic extensions are defined on the unit disk $\mbd .$ By the Riemann mapping theorem we may relax $\mbd $ to a Jordan domain $\mbx .$ Still denote by $h$ the harmonic diffeomorphism from $\mbx$ onto a bounded convex Jordan domain $\mby .$   
The Sobolev regularity for such $h$ is also interesting. Assume that $\mby$ is the unit disk, and $\mbx$ has $s$-hyperbolic growth for $s \in (0,1),$ see Definition \ref{1230_1} for the precise definition. 
Koski and Onninen recently proved in \cite[Theorem 3.2]{Koski 2018} that  $h \in W^{1,p} (\mbx, \mbd)$ for all $p < 1+s .$ In  this paper, we explore the estimate in Orlicz classes near $L^{1+s} (\mbx)$ for the derivatives of $h$ and the weighted $W^{1,1+s}$ estimate of $h .$ Our first main result is the following.

\begin{thm}\label{thm1}
Let $\X$ be a Jordan domain with $s$-hyperbolic growth for $s \in (0,1)$ and $\varphi: \partial\X\rightarrow \partial\D$ be a homeomorphism. Let  $h: \X\rightarrow \D$ be the harmonic extension of $\varphi$. Then $h\in W^{1, \Phi}(\X,\mathbb D)$ where $\Phi(t)=t^{1+s}\log^\lambda(e+t) $ for all $\lambda<-1$.
\end{thm}

Our second main result concerns weighted $W^{1,1+s}$ estimate of $h .$ 
Denote by $d(z,\partial \mbx)$ the Euclidean distance between $z$ and the boundary $\partial \mbx$ of a Jordan domain $\mbx.$ 
\begin{thm}\label{1123_1}
Let $\mbx$ be a Jordan domain with $s$-hyperbolic growth for $s \in (0,1)$ and $\varphi :\partial \mbx \ra \partial \mbd $ be a homeomorphism. 
The harmonic extension $h :\mbx \ra \mbd$ of $\varphi$ satisfies that $\int_{\mbx} |Dh (z)|^{1+s} \log^{\lambda} (e+ \frac{1}{d(z, \partial \mbx)})\, dz <\infty$ for all $\lambda <-1 .$
\end{thm}

Example \ref{20201122} shows the sharpness of the ranges on $\lambda$ in both Theorem \ref{thm1} and Theorem \ref{1123_1}. In Corollary \ref{0113_3} we obtain analogous integrability results on more general domains. Moreover Example \ref{0126_1} shows the optimality of these regularity.

Besides Poisson extensions, there are other methods to construct extensions. Let $\mby$ be a quasidisk.
Given a homeomorphism $\varphi :\mbs \ra \partial \mby,$
Koskela-Koski-Onninen\cite{KoKoOn20} found a $W^{1,p} _{\loc}(\mbr^2 ,\mbr^2)$ homeomorphic extension for all $p<2 .$
They used the bi-Lipschitz characterization of quasidisks by Rohde\cite{Rohde 2001} and an extension result of quasidisks by Tukia\cite{Tukia 1980}.
When $\mby$ is a John disk, same arguments imply a $W^{1,p} (\mbd ,\mby)$ homeomorphic extension for all $p<2 .$ If $\mby$ is a quasidisk, the complement domain $\mby^c$ is also a quasidisk. This two-sides property does not hold for John disks. Xu et al. \cite{Xu20b, GoXu20} discussed extension problems on the complement domain of a particular John disk. Assume that $\mby$ is the standard cardioid domain 
\begin{equation*}
\mby=\left\{(x, y) \in \mbr^2 : (x^2 +y^2)^2 -4x(x^2 +y^2)-4y^2<0 \right\},
\end{equation*}
and $\varphi :\mbs \ra \partial \mby$ is a homeomorphism induced by a conformal mapping from $\mbd$ onto $\mby$ via the Osgood-Caratheodory theorem. Xu constructed in \cite{Xu20b} a Sobolev homeomorphic extension of $\varphi$ on $\mbr^2 .$
From the cardioid geometry of $\mby,$ he decomposes $\mby^c$ into a sequence of pieces, and then constructs an extension on every piece.
Later with Guo, Xu in \cite{GoXu20} generalized this argument to the case that $\varphi :\mbs \ra \partial \mby$ is a homeomorphism induced by a quasiconformal mapping from $\mbd$ onto $\mby .$
Recently Koski and Onninen \cite{KoOn20} introduced another extension method. Without additional requirements, $\mby$ is just a Jordan domain. They decompose $\mby$ by hyperbolic geodesics into a sequence of pieces. Afterwards they construct an extension on each piece by Kovalev's bi-Lipschitz extension results in \cite{Kovalev 2019}. The aimed extension $\Phi :\mbr^2 \ra \mbr^2$ is the combination of extensions on pieces. The Sobolev integral of $\Phi$ has an upper bound, which is a sum with infinite terms and each term relates the Lipschitz constant on one piece.

The paper is organized as follows. In Section \ref{pre}, we present preliminaries about the hyperbolic metric, the quasihyperbolic metric and Young functions.
In Section \ref{prf main} we show the proofs of Theorem \ref{thm1} and Theorem \ref{1123_1}, and generalize the statements of both theorems to Corollary \ref{0113_3}. We construct counter-examples in Section \ref{cexm} to show the optimality of the preceding results.

\textbf{Notation.} 
For $i = 1, 2, ...$, denote $\log_{(i)}(x) = \log \left(\cdots \left(\log \left (\log(x)\right) \right)\cdots \right)$ is $i$-iterated logarithmic function and and $e_{i} = \exp  \left(\cdots \left( \exp \left ( \exp \right) \right)\cdots \right)$ is $i$-iterated exponent.

\section{Preliminaries}\label{pre}

In this section we review the hyperbolic metric, the quasihyperbolic metric, and Young functions.

The hyperbolic geometry, also called Lobachevsky-Bolyai-Gauss geometry, satisfies all of Euclid's postulates except the parallel postulate. The discovery of hyperbolic geometry is a major breakthrough in mathematics.
We can define hyperbolic geometry by axiomatic foundations. Then it is not easy to access basic and elementary results in hyperbolic geometry. We shall alternatively describe hyperbolic geometry in terms of Euclidean geometry. Hence hyperbolic geometry can be seen as a subordinate of Euclidean geometry.
Points, lines and other configurations will be defined as subsets of Euclidean geometry.

\begin{defn}
The hyperbolic metric $j_{\mbd}$ on $ \mbd $ is defined by 
\begin{equation}\label{0112_1}
j_{\mbd} (z_1 ,z_2) = \inf_{\gamma} \int_{\gamma} \frac{2 \, |dz|}{1-|z|^2} ,
\end{equation}
where the infimum is taken from all rectifibale arcs $\gamma$ joining $z_1$ and $z_2$ in $\mbd .$
\end{defn}

It is well-known that 
\begin{equation}\label{0112}
j_{\mbd} (z_1 ,z_2) = \log \frac{|1-\bar{z}_1 z_2| +|z_1 -z_2|}{|1-\bar{z}_1 z_2| -|z_1 -z_2|} .
\end{equation}
Let $\Omega$ be a simply connected domain. By the Riemann mapping theorem, we can define the hyperbolic metric on $\Omega .$ In fact, by a conformal mapping $f$ we can map $\mbd$ onto $\Omega .$ Afterwards the hyperbolic metric $j_{\Omega}$ on $\Omega$ is defined as 
\begin{equation}\label{0112_3}
j_{\Omega} (z_1 ,z_2) := j_{\mbd} (f^{-1} (z_1), f^{-1} (z_2)).
\end{equation} 
The good definition of $j_{\Omega}$ comes from the fact that $j_{\mbd}$ is invariant under M\"obius transformations.
Usually we do not have the explicit formula of the above $f .$ Therefore we cannot explicitly calculate formula for $j_{\Omega},$ like the one for $j_{\mbd}$ as in \eqref{0112}. 
Next we recall a useful substitute for $j_{\Omega} ,$ whose definition is analogous to \eqref{0112_1}.

\begin{defn}
The quasihyperbolic metric $h_{\Omega}$ on a domain $\Omega$ is defined as 
\begin{equation}\label{0113_5}
h_{\Omega} (z_1 ,z_2) = \inf_{\gamma} \int_{\gamma} \frac{1}{d(z, \pao)} \, |dz|
\end{equation}
where the infimum is taken from all rectifibale arcs $\gamma$ joining $z_1$ and $z_2$ in $\Omega ,$ and $d(z, \pao)$ is the euclidean distance between $z$ and $\pao .$
\end{defn}

In order to study the quasiconformal homogeneity, Gehring and Palka \cite{Gehring 1976} introduced quasihyperbolic metric, and then this metric shows a number of applications. For example, P. W. Jones \cite{Jones 1980} use this metric to study extension theorems of BMO functions. Domains, onto which conformal mappings of the unit disk are global H\"older continuous, are explicitly characterized by an analytic condition on this metric, see Becker-Pommerenke's work \cite{Becker 1982}. We recommend an expert survey on this metric by P. Koskela \cite{Koskela 1995}.

We show how $h_{\Omega}$ works as a substitute for $j_{\Omega}.$
A corollary of the Koebe distortion theorem (\cite[Corollary 1.4]{Pommerenke 1992}) states that a conformal mapping $g :\mbd \ra \mbc$ satisfies that
\begin{equation}\label{0112_4}
|g' (z) | \approx \frac{d(g(z) ,\partial g(\mbd))}{d(z, \partial \mbd)}
\end{equation}
for all $z \in \mbd .$ Therefore when $\Omega$ is simply connected, from \eqref{0112_1}, \eqref{0112_3} and \eqref{0113_5} the Riemann mapping theorem implies that 
\begin{equation}\label{0112_2}
h_{\Omega} (z_1 , z_2) \approx j_{\Omega} (z_1 ,z_2)
\end{equation}
for all $z_1 \in \Omega$ and $z_2 \in \Omega .$

By using the quasihyperbolic metric, the domains with hyperbolic growth are defined.  
\begin{defn}\label{1230_1} 
Let $\Omega$ be a planar domain. Fix a point $z_0 \in \Omega .$ We say that $\Omega$ satisfies $s$-hyperbolic growth for $s \in (0,1),$ if 
\begin{equation}\label{1230_2}
h_{\Omega} (z_0 ,z) \le \left(\frac{d(z_0 ,\pao)}{d(z, \pao)} \right)^{1-s}
\end{equation}
holds for all $z \in \Omega .$
\end{defn}

Recall that $h_{\Omega} (x,y) \ge \big| \log \frac{d(x,\pao)}{d(y ,\pao)} \big|$ for all $x \in \Omega$ and $y \in \Omega$, see \cite[Lemma 2.1]{Gehring 1976}. Hence the requirement \eqref{1230_2} makes sense.
On these domains, Koski and Onninen\cite{Koski 2018} have studied Sobolev homoemorphic extensions. Their work is the motivation of this paper. 
Relaxing the power of distance as in \eqref{1230_2},  
Y. Gotoh\cite{Gotoh 2000} in 2000 studied geometric properties of these generalized domain.
And their relations to quasiconformal mappings are explored by Hencl at.al in \cite{Hencl 2005} and \cite{Koskela 2001}.

Let us show examples on domains with the hyperbolic growth condition.
John disks are an important research object in geometric analysis. It is easy to check that a $c$-John disk $\Omega$ satisfies $h_{\Omega} (z_0 ,z) \le c^{-1} \log \frac{1}{d (z ,\pao)}$ for $c \ge 1.$ 
Secondly the degree of an outer-cusp will determine the power of distance in \eqref{1230_2}. For example, the domain 
\begin{equation*}
\Omega= \{(x,y): |y|\le x^{1/s},\ x \in [0,1] \} \cup \{(x,y): (x-1)^2 +y^2 \le 1\}
\end{equation*}
with $s \in (0,1)$ has $s$-hyperbolic growth condition.  

We next provide estimates related to conformal mappings onto domains with $s$-hyperbolic growth. They are useful for proofs of Theorem \ref{thm1} and Theorem \ref{1123_1}. 
Let $\Omega$ be a such domain, and $g: \mbd \ra \Omega$ be a conformal mapping. Then \eqref{0112_2}, \eqref{0112_3} and \eqref{0112} implies that 
\begin{equation*}
h_{\Omega} (g(0),g(x)) \approx j_{\Omega} (g(0),g(x)) \approx j_{\mbd} (0,x) \approx \log \frac{1}{1-|x|}
\end{equation*}
for all $x \in \mbd .$
In addition of \eqref{1230_2}, it follows that 
\begin{equation}\label{0112_5}
d(g (x), \pao) \lesssim \log^{\frac{1}{s-1}} \left( \frac{1}{1-|x|} \right) .
\end{equation}
Furthermore \eqref{0112_4} implies that 
\begin{equation}\label{0112_6}
|g' (x)| \lesssim \frac{1}{1-|x|} \log^{\frac{1}{s-1}} \left( \frac{1}{1-|x|} \right)
\end{equation}

A function $\Phi: [0, \infty)\rarrow [0, \infty)$ is a {\it  Young function} if 
it is a continuous, increasing and convex function satisfying  $\Phi(0)=0$,
$$\lim_{t\rarrow 0+}\frac{\Phi(t)}{t}=0\ \ \text{and}\  \  \lim_{t\rarrow +\infty}\frac{\Phi(t)}{t}=+\infty.$$
A  Young function $\Phi$ is said to satisfy the {\it $\Delta_2-$condition} if there is a constant $C_\Phi>0$ such that 
\begin{equation*}
\Phi(2 t) \le C_\Phi \Phi(t),\qquad \forall\   t \ge 0. 
\end{equation*}
Let $1<\alpha<\infty$. Then the function $\Phi(t)=t^\alpha\log^{\lambda}(e+t)$ with $\lambda\in \real$ and the function 
$$\Psi(t)=t^\alpha\log^{\sigma_1} (e+t) \log^{\sigma_2} _{(2)}(e_2+t)... \log^{\sigma_n} _{(n)}(e_n +t), \ \text{with} \ \ \sigma_i\in \real, i=1, \cdots, n$$ 
are Young functions satisfying the $\Delta_2 -$condition.

\section{Proofs of Main Results}\label{prf main}

In this section, we show  the proofs of Theorem \ref{thm1} and Theorem \ref{1123_1}. And we obtain analogies of these two theorems in more generalized domains. 
Before the proof of Theorem \ref{thm1}, we need the following theorem.
\begin{thm}\label{thm2}
Let $\X$ be a Jordan domain, and denote by $g: \D\rightarrow \X$ a conformal map onto $\X$. Let $\Phi$ be a Young function satisfying the $\Delta_2 -$condition. Suppose that the condition
\begin{equation}\label{20201228_7}
\sup_{w\in \partial\D}\int_{\D}\Phi\left(\frac{1}{|g'(z)||w-z|}\right)|g'(z)|^2\, dz\leq M<\infty
\end{equation}
holds. Then the harmonic extension $h: \X\rightarrow \D$ of any boundary homeomorphism $\varphi: \partial\X\rightarrow \partial\D$ lies in the Orlicz-Sobolev space $W^{1, \Phi}(\X, \mathbb C)$. 
\end{thm}

\begin{proof}[Proof of Theorem \ref{thm2}]
Let $g :\mbd \ra \mbx$ be a conformal mapping.
As in the method of Koski-Onninen \cite{Koski 2018}, without loss of generality, we may assume that $h \circ g$ is smooth up to the boundary. 
Hence 
\begin{equation}\label{20201228_5}
|(h \circ g)_z (z)| = \left|\int^{2 \pi} _{0} \frac{\psi' (e^{it})}{z-e^{it}} i e^{it} \, dt \right| \le \int^{2 \pi} _{0} \frac{|\psi' (e^{it})|}{|e^{it} -z|} \, dt .
\end{equation}
A change of variable and the estimate \eqref{20201228_5} imply that 
\begin{align}\label{20201228_6}
\int_{\mbx} \Phi (|h_{\tilde{z}} (\tilde{z})|) d \tilde{z}
= & \int_{\mbd} \Phi ( \frac{|(h \circ g)_z (z)|}{|g' (z)|}) |g' (z)|^2 dz \notag \\
\le  & \int_{\mbd} \Phi \left(\int^{2 \pi} _{0} \frac{|\psi' (e^{it})|}{|g' (z)| |z-e^{it}|} \, dt \right) |g' (z)|^2 \, dz.  \notag
\end{align} 
Note that $\int_{0}^{2\pi} |\psi'(e^{it})|\, dt=2\pi$. Using Jensen's inequality and the $\Delta_2 -$condition of $\Phi$, It follows from the Fubini theorem that
\begin{align*}
\int_{\mbx} \Phi (|h_{\tilde{z}} (\tilde{z})|) d \tilde{z}&\lesssim  \int_{\mbd} \Phi \left(\int^{2 \pi} _{0} \frac{|\psi' (e^{it})|}{2\pi |g' (z)| |z-e^{it}|} \, dt \right) |g' (z)|^2 \, dz \\
&= \int_{\mbd} \Phi \left(\int^{2 \pi} _{0} \frac{|\psi' (e^{it})|}{\int_0^{2\pi} |\psi'(e^{it})|\, dt \,|g' (z)| |z-e^{it}|} \, dt \right) |g' (z)|^2 \, dz\\
&\leq \int_{\mbd} \frac{1}{2\pi} \int_{0}^{2\pi} \Phi\left(\frac{1}{|g' (z)| |z-e^{it}|}\right) |\psi'(e^{it})|\, dt |g'(z)|^2\, dz\\
&=\frac{1}{2\pi}\int_0^{2\pi}  |\psi'(e^{it})| \left(\int_{\mbd} \Phi\left(\frac{1}{|g' (z)| |z-e^{it}|}\right)\,dz \right)\, dt\leq M<\infty.
\end{align*}
Hence under the condition \eqref{20201228_5}, the harmonic extension $h\in W^{1, \Phi}(\mbd, \mathbb C)$.
\end{proof}

\begin{proof}[Proof of Theorem \ref{thm1}]
From the estimate \eqref{0112_6},  we know that for any $z\in \mathbb D$,
\[|g'(z)|\leq \frac{C}{(1-|z|)\log^{1/(1-s)}\frac{1}{1-|z|}}.\]
Note that  $\Phi(t)=t^{p}\log^{\lambda}(e+t)$ with $p=1+s$ and $\lambda<-1$ is a Young function satisfying the $\Delta_2-$condition.  Since $2-p=1-s$,  we obtain that
\begin{eqnarray*}
\Phi\left(\frac{1}{|g'(z)||w-z|}\right)|g'(z)|^2 &=&\frac{|g'(z)|^{2-p}}{|w-z|^p}\log^\lambda\left(e+\frac{1}{|g'(z)||w-z|}\right)\\
&\leq& C\frac{\log^{\lambda}\left(e+\frac{(1-|z|)\log^{1/(1-s)} \frac{1}{1-|z|} }{|w-z|}\right)}{(1-|z|)^{2-p}|w-z|^p\log \frac{1}{1-|z|}}.
\end{eqnarray*}
Hence it is enough to prove the following quantity is finite as then apply the rotational symmetry will imply that the estimate \eqref{20201228_7} holds for all $\omega$.
 $$\int_{\D } \frac{\log^{\lambda}\left(e+\frac{(1-|z|)\log^{1/(1-s)} \frac{1}{1-|z|} }{|1-z|}\right)}{(1-|z|)^{2-p}|1-z|^p\log \frac{1}{1-|z|}}\, dz=: \int_{\D } F(z)\, dz =\int_{\D\setminus \frac12 \D} F(z)\, dz+\int_{\frac12 \D} F(z)\, dz.$$
 Notice that $F(z)$ is bounded on the set $\frac12 \D$, and hence the integral of $F(z)$ over $\frac12 \D$ is finite. Thus, it is sufficient to show that 
 $$\int_{\D\setminus \frac12 \D} F(z)\, dz<\infty.$$
Towards that, we follow the idea of Koski-Onninen \cite{Koski 2018} to divide $\D\setminus \frac12 \D$ into three pieces $S_1, S_2$ and $S_3$. More precisely, 
\begin{align*}
S_1 = & \{1+r e^{i \theta} : r \in [0,3/4],\ \theta \in [3 \pi /4 , 5 \pi /4]   \} , \\
S_2 = & \{(x,y) \in \mbd : y \in [-1/\sqrt{2} , 1/ \sqrt{2}],\ x \in [1-|y|,1]\} , \\
S_3 = & \{r e^{i \theta} : r \in [1/2 ,1],\ \theta \in [\pi /4 , 7 \pi /4]\} .
\end{align*}

  On the set $S^1$, we have the estimate $1-|z|\approx |1-z|$. Hence we may apply polar coordinates around the point $z=1$ to find that
 \[\int_{S_1}F(z)\, dz\lesssim \int_{S_1} \frac{\log^\lambda\left(e+\log^{1/(1-s)}\frac{1}{|1-z|}\right)}{|1-z|^2\log\frac{1}{|1-z|}}\, dz \lesssim \int_{3\pi/4}^{5\pi/4}\int_{0}^{3/4} \frac{\log^\lambda(\log(1/r))}{r\log(1/r)}\, dr\, d\theta<\infty.\]

On the set $S_3$, we have that $|1-z|$ is bounded away from zero. Hence
$\frac{1}{|1-z|^p}$ in the expression $F$ has a positive constant upper bound. Notice that
 $$\log^{\lambda}\left(e+\frac{(1-|z|)\log^{1/(1-s)} \frac{1}{1-|z|} }{|1-z|}\right) \le \log^{\lambda} (e) =1.$$
  We  have the estimate that
$$\int_{S_3}F(z)\, dz\lesssim\int_{S_3} \frac{1}{(1-|z|)^{2-p}\log \frac{1}{1-|z|}}\, dz\lesssim \int_{\pi/4}^{7\pi/4}\int_{1/2}^{1} \frac{r}{(1-r)^{2-p}}\, dr\, d\theta<\infty,$$
where the second inequality is from an estimate $\log^{-1} (\frac{1}{1-|z|}) \lesssim 1 .$

On the set $S_2$, since on each $R_\theta$, $|1-z|$ is comparable to the size of the angle $\theta$ and $1-|z|<|1-z|$, we obtain the following estimate
\[F(z)\lesssim \frac{\log^\lambda\left(e+\frac{(1-|z|)\log^{1/(1-s)}\frac{1}{|\theta|}}{|\theta|}\right)}{(1-|z|)^{2-p} |\theta|^p \log \frac{1}{|\theta|}}, \ \ z\in R_\theta.\]
By the property that $1-\rho(\theta)$ is comparable to $|\theta|$, we estimate that
\begin{eqnarray*}
\int_{S_2} F(z)\, dz&\lesssim& \int_{-\pi/4}^{\pi/4}\frac{1}{|\theta|^p\log \frac{1}{|\theta|}}\int_{\rho(\theta)}^{1}\frac{\log^{\lambda}\left(e+\frac{(1-r)\log^{1/(1-s)}\frac{1}{|\theta|}}{|\theta|}\right)}{(1-r)^{2-p}}\, dr\, d\theta\\
&\approx&  \int_{-\pi/4}^{\pi/4}\frac{1}{|\theta|^p\log \frac{1}{|\theta|}} \frac{\log^{\lambda}\left(e+\frac{(1-\rho(\theta))\log^{1/(1-s)}\frac{1}{|\theta|}}{|\theta|}\right)}{(1-\rho(\theta))^{1-p}}\, d\theta\\
&\approx&  \int_{-\pi/4}^{\pi/4} \frac{\log^\lambda(e+\log^{1/(1-s)} \frac{1}{|\theta|})}{|\theta|\log \frac{1}{|\theta|}}\,d\theta\\
&\lesssim&  \int_{-\pi/4}^{\pi/4} \frac{\log^\lambda(\log(1/|\theta|))}{|\theta|\log(1/|\theta|)}\, d\theta<\infty.
\end{eqnarray*}
Summing up the integrations over the sets $S_1$, $S_2$ and $S_3$, the proof is finished.
\end{proof}

Towards the  proof of Theorem \ref{1123_1}, we give the following lemma first.

\begin{lem}\label{1123}
Let $\mbx$ be a Jordan domain with $s$-hyperbolic growth for $s \in (0,1)$. For any conformal mapping $g :\mbd \ra \mbx,$ we have that 
\begin{equation*}
\int_{\mbd} \frac{|g'(z)|^{1-s}}{|1-z|^{1+s}} \log^{\lambda} \log (\frac{1}{1-|z|} ) \, dz <\infty
\end{equation*}
for all $\lambda <-1 .$
\end{lem}

\begin{proof}[Proof of Lemma \ref{1123}]
The inequality as in \eqref{0112_6} states that $|g' (z)| \lesssim \frac{1}{1-|z|} \log^{\frac{1}{s-1}} (\frac{1}{1-|z|})$ for $z \in \mbd .$ Hence replacing $|g' (z)|$ by this upper bound implies that 
\begin{equation*}
\frac{|g'(z)|^{1-s}}{|1-z|^{1+s}} \log^{\lambda} \log (\frac{1}{1-|z|} )
\lesssim  \frac{1}{|1-z|^{1+s} (1-|z|)^{1-s}} \log^{-1} (\frac{1}{1-|z|}) \log^{\lambda} \log (\frac{1}{1-|z|}). 
\end{equation*}
To complete the proof , it is then enough to check the following integral is finite
\begin{equation*}
\int_{\mbd \setminus \frac{1}{2} \mbd} \frac{1}{|1-z|^{1+s} (1-|z|)^{1-s}} \log^{-1} (\frac{1}{1-|z|}) \log^{\lambda} \log (\frac{1}{1-|z|}) \, dz =:
\int_{\mbd \setminus \frac{1}{2} \mbd} G (z) \, dz .
\end{equation*}
In sequel we let $\lambda <-1 .$
As in Koski-Onninen's work \cite{Koski 2018}, we divide $\mbd \setminus \frac{1}{2} \mbd$ into subsets $S_1 ,\ S_2$ and $S_3 .$
Here 
\begin{align*}
S_1 = & \{1+r e^{i \theta} : r \in [0,3/4],\ \theta \in [3 \pi /4 , 5 \pi /4]   \} , \\
S_2 = & \{(x,y) \in \mbd : y \in [-1/\sqrt{2} , 1/ \sqrt{2}],\ x \in [1-|y|,1]\} , \\
S_3 = & \{r e^{i \theta} : r \in [1/2 ,1],\ \theta \in [\pi /4 , 7 \pi /4]\} .
\end{align*}
Afterwards it suffices to check that $\int_{S_1} G <\infty$, $\int_{S_2} G <\infty $ and $\int_{S_3} G <\infty  .$

Notice that $|1-z| \approx 1-|z|$ for $z \in S_1 .$ Hence the polar coordinates around the point $z=1$ implies that  
\begin{align*}
\int_{S_1} G (z) \, dz \approx & 
\int_{S_1} \frac{1}{|1-z|^{2}} \log^{-1} (\frac{1}{|1-z|}) \log^{\lambda} \log (\frac{1}{|1-z|}) \, dz \notag \\
=& \int_{3 \pi/4} ^{5 \pi /4} \int_{0} ^{3/4} \frac{1}{r} \log^{-1} (\frac{1}{r}) \log^{\lambda} \log (\frac{1}{r}) \, dr \, d\theta <\infty .
\end{align*}

On the set $S_3 ,$ the modulus $|1-z|$ is bounded away from zero. 
Therefore by bounding the term $|1-z|^{-1-s}$ in the integrand $G$ and applying polar coordinates around the origin, we have that 
\begin{align*}
\int_{S_3} G (z) \, dz \lesssim &
\int_{S_3} \frac{1}{ (1-|z|)^{1-s}} \log^{-1} (\frac{1}{1-|z|}) \log^{\lambda} \log (\frac{1}{1-|z|}) \, dz \notag \\
= & \int_{\pi /4} ^{7 \pi /4} \int_{1/2} ^1  \frac{r}{ (1-r)^{1-s}} \log^{-1} (\frac{1}{1-r}) \log^{\lambda} \log (\frac{1}{1-r}) \, dr \, d \theta <\infty .
\end{align*}

Finally, we check that $\int_{S_2} G <\infty .$
Applying the triangle inequality $1-|z| \le |1-z|,$ by $\log^{-1} \frac{1}{|1-z|}$ we can replace the term $\log^{-1} \frac{1}{1-|z|}$ in the integrand $G .$ Afterwards we use the polar coordinates $(r, \theta)$ around the origin.
The modulus $|1-z|$ is comparable to the size of the angle $\theta .$
For each angle $\theta,$ the modulus $r=|z|$ ranges from $\rho (\theta)$ to $1 .$ And $1-\rho(\theta) = \frac{\sin (\pi /4 +\theta) -\sin(\pi /4)}{\sin (\pi /4 +\theta}$ is comparable to $|\theta|.$
By these estimates, we obtain that 
\begin{align*}
\int_{S_2} G (z) \, dz \lesssim &
\int_{S_2} \frac{1}{(1-|z|)^{1-s}} \frac{1}{|1-z|^{1+s} } \log^{-1} (\frac{1}{|1-z|}) \log^{\lambda} \log (\frac{1}{|1-z|}) \, dz \notag \\
\lesssim & \int_{-\pi/4} ^{\pi /4} \frac{1}{|\theta|^{1+s} } \log^{-1} (\frac{1}{|\theta|}) \log^{\lambda} \log (\frac{1}{|\theta|}) \, d \theta \int_{\rho (\theta)} ^1 \frac{1}{(1-r)^{1-s}} \, dr  \notag \\
\approx & \int_{-\pi/4} ^{\pi /4} \frac{1}{|\theta|} \log^{-1} (\frac{1}{|\theta|}) \log^{\lambda} \log (\frac{1}{|\theta|}) \, d \theta <\infty  .
\end{align*}
\end{proof}

\begin{proof}[Proof of Theorem \ref{1123_1}]
Take a conformal mapping $g : \mbd \ra \mbx .$ Then \eqref{0112_5} shows that $d(g(z) , \partial \mbx) \lesssim \log^{\frac{1}{s-1}} (\frac{1}{1-|z|}) $ for all $z \in \mbd.$
Together with a change of variables $w=g(z)$, we estimate that 
\begin{align}\label{20201227}
\int_{\mbx} |h' (w)|^{1+s} \log^{\lambda} (e+\frac{1}{d(w, \partial \mbx)}) dw = & \int_{\mbd} |(h \circ g)' (z)|^{1+s} |g' (z)|^{1-s}  \log^{\lambda} (e+ \frac{1}{d(g(z),\partial \mbx)}) dz \notag\\
\lesssim & \int_{\mbd} |(h \circ g)' (z)|^{1+s} |g' (z)|^{1-s} \log^{\lambda} \log (\frac{1}{1-|z|}) dz .
\end{align}
Next we look for an upper bound for the rightest integral of \eqref{20201227}.

Notice that $h \circ g$ is a harmonic mapping, which can be expressed as the Poisson extension $\frac{1}{2 \pi}\int_{\mbs} \frac{1-|z|^2}{|z-\xi|^2} \psi (\xi) \, |d \xi| .$ Here $\psi = \varphi \circ g|_{\mbs} $ is a self-homeomorphism on $\mbs .$ Basically $\psi$ can be seen as a self-homeomorphism (then a strictly monotonic function) on the interval $[0,2 \pi].$ Hence the derivative $\psi'$ exists a.e. on $\mbs$ and 
\begin{equation}\label{0112_7}
\int^{2 \pi} _{0} |\psi' (e^{it})| \, dt <\infty .
\end{equation}   
Furthermore by differentiating the Poisson integral of $h \circ g ,$ we obtain a pointwise estimate on the derivative of $h \circ g$
\begin{equation*}
|(h \circ g)' (z) |= \big| \int^{2 \pi} _{0} \frac{\psi' (e^{it})}{z- e^{it}} i e^{it} \, dt \big| \le \int^{2 \pi} _{0} \frac{|\psi' (e^{it})|}{|z- e^{it} |} \, dt 
\end{equation*}
for all $z \in \mbd .$
By this estimate and Minkowski's inequality, we estimate
\begin{align}\label{20201227_1}
& \int_{\mbd} |(h \circ g)' (z)|^{1+s} |g' (z)|^{1-s} \log^{\lambda} \log (\frac{1}{1-|z|}) dz  \notag\\
\le & \int_{\mbd} \left( \int^{2 \pi} _{0} \frac{|\psi' (e^{it})|}{|z- e^{it}|} \, dt \right)^{1+s} |g' (z)|^{1-s} \log^{\lambda} \log (\frac{1}{1-|z|}) \, dz  \notag \\
\le & \left( \int^{2 \pi} _{0} |\psi' (e^{it})| \left( \int_{\mbd} \frac{|g' (z)|^{1-s}}{|z-e^{it}|^{1+s}} \log^{\lambda} \log (\frac{1}{1-|z|})\, dz \right)^{1/(1+s)} \, dt \right)^{1+s} .
\end{align}
Composition of \eqref{20201227} with \eqref{20201227_1} implies that 
\begin{align}\label{20201227_2}
& \int_{\mbx} |h' (w)|^{1+s} \log^{\lambda} (e+\frac{1}{d(w, \partial \mbx)})\, dw  \notag \\
\lesssim & \left( \int^{2 \pi} _{0} |\psi' (e^{it})| \left( \int_{\mbd} \frac{|g' (z)|^{1-s}}{|z-e^{it}|^{1+s}} \log^{\lambda} \log (\frac{1}{1-|z|})\, dz \right)^{1/(1+s)} \, dt \right)^{1+s} .
\end{align}
By the rotation invariance and Lemma \ref{1123}, we have that 
\begin{equation*}
\int_{\mbd} \frac{|g' (z)|^{1-s}}{|z-e^{it}|^{1+s}} \log^{\lambda} \log (\frac{1}{1-|z|})\, dz 
=  \int_{\mbd} \frac{|g'(z)|^{1-s}}{|1-z|^{1+s}} \log^{\lambda} \log (\frac{1}{1-|z|} ) \, dz <\infty 
\end{equation*}
for all $t \in [0, 2 \pi] .$
Therefore in addition of \eqref{0112_7}, we derive from \eqref{20201227_2} that 
\begin{equation}
\int_{\mbx} |h' (w)|^{1+s} \log^{\lambda} (e+\frac{1}{d(w, \partial \mbx)})\, dw 
\lesssim  \left( \int^{2 \pi} _{0} |\psi' (e^{it})| \, dt \right)^{1+s} <\infty .
\end{equation}
The proof is complete.
\end{proof}


In Theorem \ref{thm1} and Theorem \ref{1123_1}, we obtain the regularity of harmonic extension on domains satisfying the condition \eqref{1230_2}. Now we relax these domains via replacing the upper bound $(\frac{d(z_0 ,\pao)}{d(z, \pao)} )^{1-s}$ in \eqref{1230_2} by $(\frac{1}{d(x, \partial \mbx)} )^{1-s} \log^{\sigma} (e+\frac{1}{d(x, \partial \mbx)}) .$ 
We hope to discuss analogous regularity on these domains.

\begin{cor}\label{20201228}
Let $\mbx$ be a Jordan domain satisfying that $h_{\mbx} (x_0 ,x) \le (\frac{1}{d(x, \partial \mbx)} )^{1-s} \log^{\sigma} (e+\frac{1}{d(x, \partial \mbx)}) $ for $\sigma \ge -1 .$ Let $\varphi :\partial \mbx \ra \mbs$ be a homeomorphism and $h : \mbx \ra \mbd$ be its harmonic extension. Suppose that $\Phi(t)=t^{1+s}\log^\lambda(e+t)$. Then
\begin{equation}\label{eq-314-1}
\int_{\mathbb X} \Phi(|Dh(z)|)\, dz<\infty
\end{equation}
 and 
 \begin{equation}\label{eq-314-2}
 \int_{\mbx} |Dh (z)|^{1+s} \log^{\lambda} (e+ \frac{1}{d(z, \partial \mbx)})\, dz<\infty
 \end{equation}
 whenever $\lambda <-1-\sigma .$
\end{cor}

\begin{proof}[Proof of Corollary \ref{20201228}]
The proof is analogous to that for Theorem \ref{thm1} and Theorem \ref{1123_1}, hence we only sketch the process. 
Let $g :\mbd \ra \mbx$ be a conformal mapping. 
Under the assumption $h_{\mbx} (x_0 ,x) \le \Psi (\frac{1}{d(x, \partial \mbx)})$ with $\Psi(t)=t^{1-s}\log^\sigma(e+t)$, same arguments as for \eqref{0112_5} and \eqref{0112_5} imply that 
\begin{align}
& d(g(z) ,\partial \mbx) \lesssim \log^{\frac{1}{s-1}} (\frac{1}{1-|z|}) \log^{\frac{\sigma}{1-s}} \log (\frac{1}{1-|z|}) \label{0113_1}\\
& |g' (z)| \approx \frac{d(g(z) ,\partial \mbx)}{1-|z|} \lesssim \frac{1}{1-|z|} \log^{\frac{1}{s-1}} (\frac{1}{1-|z|}) \log^{\frac{\sigma}{1-s}} \log (\frac{1}{1-|z|}) \label{0113_2}
\end{align}
for all $z \in \mbd .$
Here we used the fact  that $\Psi^{-1} (t) \approx t^{\frac{1}{1-s}} \log^{\frac{ \sigma}{s-1}} (t) .$
The estimate \eqref{eq-314-1} is obtained by applying   \eqref{0113_2} into the proof of Theorem \ref{thm1}. 

Applying \eqref{0113_1} and \eqref{0113_2} in order, from a change of variable $w=g(z)$ and Minkowski's inequality we obtain that 

\begin{align}\label{20201228_4}
& \int_{\mbx} |h' (w)|^{1+s} \log^{\lambda} (e+\frac{1}{d(w, \partial \mbx)}) dw \notag \\
\lesssim  & \int_{\mbd} \frac{|g' (z)|^{1-s}}{|1-z|^{1+s}} \log^{\lambda} \log (\frac{1}{1-|z|})\, dz    \notag \\
\lesssim & \int_{\mbd} \frac{1}{|1-z|^{1+s} (1-|z|)^{1-s}} \log^{-1} (\frac{1}{1-|z|}) \log^{\sigma+\lambda} \log (\frac{1}{1-|z|})\, dz .
\end{align}
Analogously to the proof of Lemma \ref{1123}, the lowest integral in \eqref{20201228_4} is finite whenever $\lambda <-1-\sigma,$ which gives the estimate \eqref{eq-314-2}.  The proof is complete.
\end{proof}

Here we give a more generalized variant of Corollary \ref{20201228}.
Let $\sigma = (\sigma_1 ,...,\sigma_n)$ be a vector with $\sigma_i \in \mbr$ for all $i=1,...,n .$ For $s \in \mbr$, denote 
\begin{equation}\label{0126}
\Psi_{1-s , \sigma} (t) = t^{1-s} \log^{\sigma_1} (e+t) \log^{\sigma_2} _{(2)}(e_2+t)... \log^{\sigma_n} _{(n)}(e_n +t) 
\end{equation}

\begin{cor}\label{0113_3}
Let $\mbx$ be a Jordan domain satisfying that $h_{\mbx} (x_0 ,x) \le \Psi_{1-s ,\sigma} (\frac{1}{d(x, \partial \mbx)}) $ with $s \in (0,1)$ and $\sigma_i \ge -1$ for all $i =1,...,n .$ Let $\varphi :\partial \mbx \ra \mbs$ be a homeomorphism and $h : \mbx \ra \mbd$ be its harmonic extension. Then for $\lambda=(\lambda_1 ,...,\lambda_n)$, we have 
$$\int_{\mbx} \Psi_{1+s, \lambda}(|Dh (z)|)\, dz<\infty$$
and
\begin{equation*}
\int_{\mbx} |Dh (z)|^{1+s} \Psi_{0,\lambda} \left(  \frac{1}{d(z, \partial \mbx)}\right) \, dz <\infty
\end{equation*}
whenever $\lambda_n<-1-\sigma_n $ and $\lambda_i \le -1-\sigma_i $ for all $i=1,...,n-1 .$
\end{cor}

The proof of Corollary \ref{0113_3} is completely similar to that of Corollary \ref{20201228}. The only thing to be noticed is that for
$\Psi_{1-s , \sigma}$ as in \eqref{0126} we estimate its inverse as
\begin{align*}
\Psi^{-1} _{1-s , \sigma} (t) \approx & t^{\frac{1}{1-s}} \log^{\frac{\sigma_1}{s-1}} (e+t) \log^{\frac{\sigma_2}{s-1}}_{(2)} (e_2 +t)... \log^{\frac{\sigma_n}{s-1}}_{(n)} (e_n +t) \\
=& \Psi_{\frac{1}{s-1} , \frac{\sigma}{s-1}} (t).
\end{align*}

\section{Counter-examples related to optimal regularities}\label{cexm}

In this section, we provide an example to show the sharpness of ranges on $\lambda$ in Theorem \ref{thm1} and in Theorem \ref{1123_1}, and also an example for those in Corollary \ref{0113_3}.

\begin{example}\label{20201122}

For any $s \in (0,1),$ there is a Jordan domain $\X$ with $s$-hyperbolic growth and a homeomorphism $\varphi : \partial \X \ra \partial \mbd$ which do not admit a homeomorphic extension $h :\mbx \ra \mbd$ with $\int_{\mbx} |Dh|^{1+s} \log^{-1} (e+|Dh|)<\infty $ or $\int_{\mbx} |Dh|^{1+s} \log^{-1} (e+ \frac{1}{d(z ,\partial \mbx)} ) <\infty .$ 
\end{example}

\begin{proof}[Construction of Example \ref{20201122}]

Our construction is analogous to that by Koski-Onninen \cite{Koski 2018}.
Start from a graph $\{(x,|x|^s) : x \in [-1,1]\} .$ 
Then complete this graph to obtain a smooth Jordan curve except the above cusp point. Denote by $\mbx$ the bounded Jordan domain enclosed by the preceding curve. The above description on $\partial \mbx$ determines that $\mbx$ has $s$-hyperbolic growth. In sequel we will only care $\varphi$ and $h$ in a neighbourhood of the cusp point. 

\begin{figure}[htbp]
\centering \includegraphics[width=\textwidth]{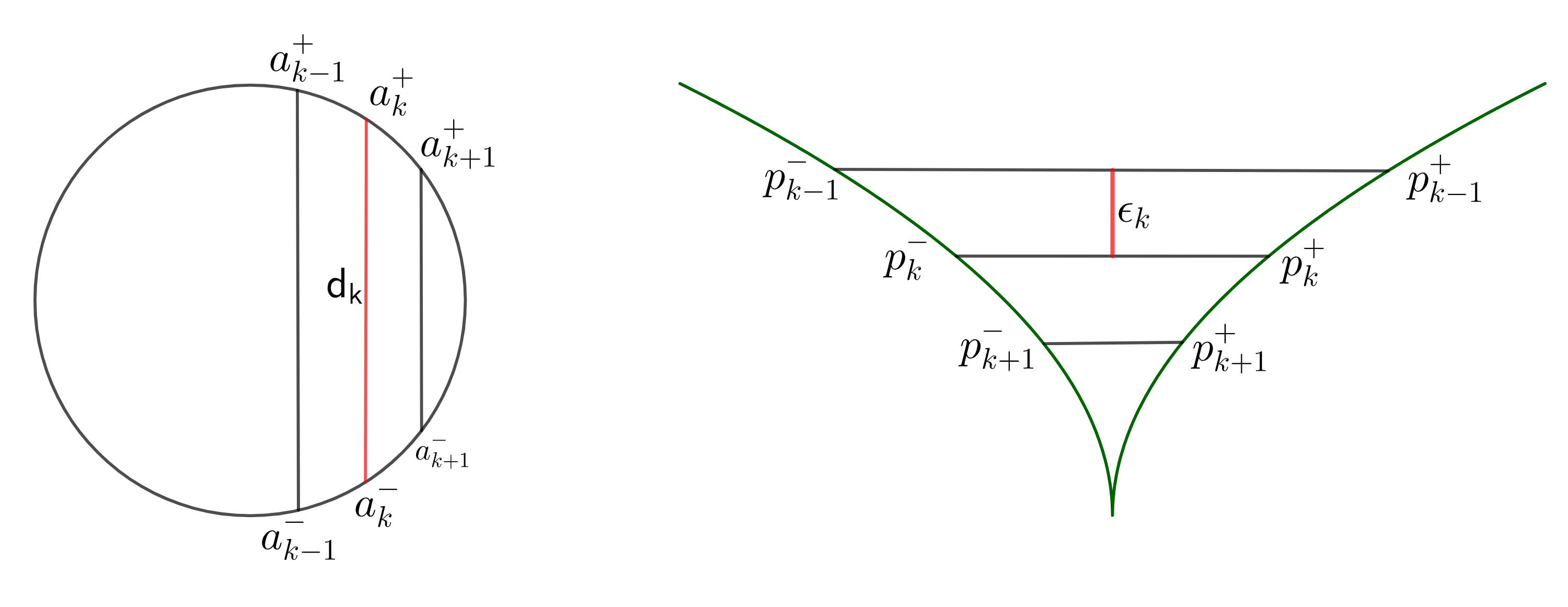}
\caption{The construction}
\label{fig1}
\end{figure}

Divide the cusp graph by a sequence of line segments parallel to the $x$-axis. Denote by $p^- _k$ the left intersection point between the $k$-line segment and the graph, and $p^+ _k$ the right one, see Figure \ref{fig1}. We require that the $y$-coordinates of $p^+ _1$ is $ \sum_{k=2} k^{-2} .$ Let $\epsilon_k = k^{-2}$ be the difference between the $y$-coordinates of $p^+ _{k-1}$ and $p^+ _{k}$ for $k=2,3,...$ 
We afterwards divide the unit disk $\mbd \cap \{(x,y): x \ge 0\}.$ To do this we use a sequence of line segments parallel to the $y$-axis. Denote by $a^+ _k$ the upper intersection point between the $k$-line segment and $\partial \mbd ,$ and $a^- _k$ the lower one. Let $d_k = \log^{\frac{-1}{1+s}} \log (1+k)$ be the length of line segment between $ a^- _k$ and $a^+ _k .$ 

By the constant speed, we define $\varphi$ to map the arc of cusp graph between $p^+ _k$ and $p^+ _{k-1}$ ($p^- _k$ and $p^- _{k-1}$) onto the arc of unit circle between $a^+ _k$ and $a^+ _{k-1} $ ($a^- _k$ and $a^- _{k-1} $). 
Assume that $h$ is a homeomorphism extension of $\varphi $ with Sobolev regularity. 
Let $S_k$ be the subset of the original cusp domain between line segments $p^- _{k-1} p^+ _{k-1}$ and $p^- _{k} p^+ _{k} .$
By Fubini's theorem and the ACL property of Sobolev functions we have that 
\begin{equation}\label{1122}
\int_{S_k} |Dh| 
\gtrsim \int_{\sum_{j \ge k+1} \epsilon_j} ^{\sum_{j \ge k}\epsilon_j }  \int^{y^{1/s}} _{-y^{1/s}} \big| \frac{\partial h (x,y)}{\partial x} \big| \, dx \, dy
\gtrsim \int_{\sum_{j \ge k+1} \epsilon_j} ^{\sum_{j \ge k} \epsilon_j}  d_k   \, dy = \epsilon_k d_k .
\end{equation}
Here we only care about the sufficiently large $k .$
Note that $|S_k| \approx \epsilon_k (\sum^{\infty} _{j = k}\epsilon_j)^{1/s} \approx k^{-2-\frac{1}{s}} .$
Hence \eqref{1122} implies that 
\begin{equation}\label{0113_4}
\displaystyle  \mint{-}_{S_k} |Dh| \gtrsim k^{\frac{1}{s}} \log^{ \frac{-1}{1+s}} \log (k) .
\end{equation}
The function $t^{1+s} \log^{-1} (e+t)$ is convex and increasing whenever $t \gg 1.$
Therefore Jensen's inequality and the estimate \eqref{0113_4} in order imply that 
\begin{align*}
\int_{S_k} |Dh|^{1+s} \log^{-1} (e+|Dh|) 
\ge & |S_k | \left(\displaystyle  \mint{-}_{S_k} |Dh| \right)^{1+s} \log^{-1} \left(e+\displaystyle  \mint{-}_{S_k} |Dh| \right) \\
\gtrsim &|S_k| \left( k^{\frac{1}{s}} \log^{ \frac{-1}{1+s}} \log (k) \right)^{1+s} \log^{-1} (k^{\frac{1}{s}} \log^{\frac{-1}{1+s}} \log (k)) \\
\approx & \frac{1}{k \log (k) \log \log (k)} .
\end{align*}
Afterwards 
\begin{align*}
\int_{\cup^{\infty} _{k = 1} S_k} |Dh|^{1+s} \log^{-1} (e+|Dh|) 
= & \sum^{\infty} _{k = 1} \int_{S_k} |Dh|^{1+s} \log^{-1} (e+|Dh|) \\
\gtrsim & \sum^{\infty} _{k = 1} \frac{1}{k \log (k) \log \log (k)} =\infty .
\end{align*}

We next check that $\int_{\mbx} |Dh|^{1+s} \log^{-1} \left(e+ \frac{1}{d(x ,\partial \mbx)} \right)=\infty .$ 
By H\"older's inequality we have 
\begin{align*}
\int_{S_k} |Dh| 
= & \int_{S_k} |Dh(z)| \log^{\frac{-1}{1+s}} (e+\frac{1}{d(z ,\partial \mbx)}) \log^{\frac{1}{1+s}} (e+\frac{1}{d(z ,\partial \mbx)}) \, dz\\
\le &\left(\int_{S_k} |Dh (z)|^{1+s} \log^{-1} (e+ \frac{1}{d(z ,\partial \mbx)} ) \, dz \right)^{\frac{1}{1+s}} 
\left(\int_{S_k} \log^{\frac{1}{s}} (e+ \frac{1}{d(z ,\partial \mbx)} ) \, dz\right)^{\frac{s}{1+s}}
\end{align*}
Therefore 
\begin{equation}\label{1122_1}
\int_{S_k} |Dh|^{1+s} \log^{-1} (e+ \frac{1}{d(z ,\partial \mbx)} )
\ge \left(\int_{S_k} |Dh| \right)^{1+s}
\left( \int_{S_k} \log^{\frac{1}{s}} (e+ \frac{1}{d(x ,\partial \mbx)} ) \right)^{-s}
\end{equation}
Notice that 
\begin{align*}
\int_{S_k} \log^{\frac{1}{s}} (e+ \frac{1}{d(z ,\partial \mbx)} ) \, dz 
\approx & \int_{\sum_{j \ge k+1} \epsilon_j} ^{\sum_{j \ge k} \epsilon_j} 
\int_{0} ^{y^{1/s}} \log^{\frac{1}{s}} \left(e+ \frac{1}{y^{1/s}-x} \right)  \, dx \, dy \\
\approx & \int_{\sum_{j \ge k+1} \epsilon_j} ^{\sum_{j \ge k} \epsilon_j}  y^{1/s} \log^{1/s} (e+ y^{-1/s}) \, dy \\
\approx & \frac{1}{k^{2+\frac{1}{s}}} \log^{1/s} ( k) .
\end{align*} 
Therefore in addition of \eqref{1122}, we obtain from \eqref{1122_1} that 
\begin{equation*}
\int_{S_k} |Dh|^{1+s} \log^{-1} (e+ \frac{1}{d(z ,\partial \mbx)} ) \, dz
\ge (\epsilon_k d_k)^{1+s} \left(\frac{1}{k^{2+\frac{1}{s}}} \log^{1/s} ( k) \right)^{-s} = \frac{1}{k \log (k) \log \log (k)} .
\end{equation*}
Afterwards 
\begin{align*}
\int_{\cup^{\infty} _{k = 1} S_k} |Dh|^{1+s} \log^{-1} (e+ \frac{1}{d(z ,\partial \mbx)} )
= & \sum^{\infty} _{k = 1} \int_{S_k} |Dh|^{1+s} \log^{-1} (e+ \frac{1}{d(z ,\partial \mbx)} )\\
\gtrsim & \sum^{\infty} _{k = 1} \frac{1}{k \log (k) \log \log (k)} =\infty .
\end{align*}
\end{proof}


The following example is to show the sharpness of $\lambda_i$ in Corollary \ref{0113_3}.
\begin{example}\label{0126_1}
Let $\Psi_{1-s , \sigma}$ be as in \eqref{0126}.
We construct a Jordan domain $\mbx$ satisfying that $h_{\mbx} (x_0 ,x) \le \Psi_{1-s ,\sigma} (\frac{1}{d(x, \partial \mbx)}) $ with $s \in (0,1) .$
There is a homeomorphism $\varphi : \partial \mbx \ra \mbs,$ which does not admit a homeomorphic extension $h: \mbx \ra \mbd$ with $\int_{\mbx} \Psi_{1+s ,\lambda} (|Dh|)<\infty$ or $\int_{\mbx} |Dh|^{1+s} (z) \Psi_{0, \lambda} (\frac{1}{d(z, \partial \mbx)})\, dz <\infty$ when $\lambda=(\lambda_1 ,...,\lambda_n)$ with $\sigma_i +\lambda_i = -1$ for any $i=1,...,n .$  
\end{example}

\begin{proof}
The construction of $\mbx$ and $\varphi$ in Example \ref{0126_1} is analogous to that in Example \ref{20201122}. In sequel we sketch the process and estimates, and leave detailed arguments to interested reader. We follow the notation in Example \ref{20201122}.   
By the cusp graph $\{(x,\Psi_{-s ,\sigma} (1/|x|)): x \in [-1 ,1]\} ,$ we define a Jordan domain $\mbx .$ We calculate that $\mbx$ satisfies $h_{\mbx} (z_0 ,z) \le \Psi_{1-s, \sigma} (1/d (z, \partial \mbx)) .$
Let 
\begin{equation*}
\epsilon_k = \frac{1}{k^2} \mbox{ and } d_k = \log^{\frac{-1}{1+s}} _{(n+1)} (e_{n+1}+k) .
\end{equation*}
Then we piecewise define $\varphi :\partial \mbx \ra \mbs $ as in Example \ref{20201122}. Assume that $h :\mbx \ra \mbd$ is a Sobolev homeomorphic extension of $\varphi .$ On a piece $S_k \subset \mbx$ we estimate that 
\begin{equation}\label{0125}
\int_{S_k} |Dh| \gtrsim \epsilon_k d_k 
\mbox{ and }
|S_k| \approx  \epsilon_k \Psi_{-\frac{1}{s}, - \frac{\sigma}{s}} (k).
\end{equation}
Therefore
\begin{align*}
\int_{S_k} \Psi_{1+s ,\lambda} (|Dh|) \ge & |S_k|  \Psi_{1+s ,\lambda} (\displaystyle  \mint{-}_{S_k} |Dh|) \\
\gtrsim & \frac{1}{k} \log^{\sigma_1 +\lambda_1} (k) ... \log^{\sigma_n +\lambda_n} _{(n)} (k) \log^{-1} _{(n+1)} (k)
\end{align*}
Finally $\int_{\mbx} \Psi_{1+s ,\lambda} (|Dh|) \gtrsim \sum^{\infty} _{k=1} \frac{1}{k} \log^{\sigma_1 +\lambda_1} (k) ... \log^{\sigma_n +\lambda_n} _{(n)} (k) \log^{-1} _{(n+1)} (k) =\infty$ whenever $\sigma_i +\lambda_i = -1$ for any $i=1,...,n .$

As in \eqref{1122_1} we estimate that 
\begin{equation}\label{0125_1}
\int_{S_k} |Dh|^{1+s} \Psi_{0, \lambda} (\frac{1}{d(z, \partial \mbx)})\, dz \ge \left( \int_{S_k} |Dh|\right)^{1+s} \left( \int_{S_k} \Psi_{0, -\frac{\lambda}{s}} (\frac{1}{d(z, \partial \mbx)} ) \, dz \right)^{-s} .
\end{equation}
Moreover we calculate that 
$\int_{S_k} \Psi_{0, -\frac{\lambda}{s}} (\frac{1}{d(z, \partial \mbx)} ) \, dz 
\approx  \Psi_{-2-\frac{1}{s} , - \frac{\sigma +\lambda}{s}} (k) .$
Here $\frac{\sigma +\lambda}{s} = (\frac{\sigma_1 +\lambda_1}{s},...,\frac{\sigma_n +\lambda_n}{s}) .$
Together with the estimate $\int_{S_k} |Dh|$ as in \eqref{0125}, we derive from \eqref{0125_1} that 
\begin{equation*}
\int_{S_k} |Dh|^{1+s} (z) \Psi_{0, \lambda} (\frac{1}{d(z, \partial \mbx)})\, dz
\gtrsim \frac{1}{k} \log^{\sigma_1 +\lambda_1} (k) ... \log^{\sigma_n +\lambda_n} _{(n)} (k) \log^{-1} _{(n+1)} (k) .
\end{equation*}
Finally $\int_{\mbx} |Dh|^{1+s} (z) \Psi_{0, \lambda} (\frac{1}{d(z, \partial \mbx)})\, dz \gtrsim \sum^{\infty} _{k=1} \frac{1}{k} \log^{\sigma_1 +\lambda_1} (k) ... \log^{\sigma_n +\lambda_n} _{(n)} (k) \log^{-1} _{(n+1)} (k) =\infty$ whenever $\sigma_i +\lambda_i = -1$ for any $i=1,...,n .$ 
\end{proof}

\section*{Acknowledgements}
Haiqing Xu is funded by the postdoctor foundation at Shandong University (No. 10000072110302), and partially supported by the Qilu funding of Shandong University (No. 62550089963197).

\bigskip
\bibliographystyle{amsplain}

\medskip

\noindent Zhuang Wang,

\noindent
MOE-LCSM, School of Mathematics and Statistics, Hunan Normal University, Changsha, Hunan 410081, P. R. China.

\noindent{\it E-mail address}:  \texttt{zhuang.z.wang@foxmail.com}, \texttt{zwang@hunnu.edu.cn}
\bigskip

\noindent Haiqing Xu,

\noindent
Research Center for Mathematics and Interdisciplinary
Sciences, Shandong University, Qingdao, Shandong 266237,
P. R. China.

\noindent{\it E-mail address}:  \texttt{hqxu@mail.ustc.edu.cn}
\end{document}